\documentclass[reqno]{amsart}

\usepackage{amsfonts, amssymb,amsmath,amsthm}
\usepackage{mathtools}
\usepackage{float}
\usepackage{hyperref}
\hypersetup{colorlinks = true}

\usepackage{color}

\renewcommand{\ss}{\subsection}

\renewcommand{\d}{\partial}
\newcommand{\pa}{\partial}

\newcommand{\la}{\langle}
\newcommand{\ra}{\rangle}
\newcommand{\ot}{\otimes}

\newcommand{\h}{\hat}
\newcommand{\wt}{\widetilde}

\newcommand{\wb}{\overline}

\newcommand{\z}{\text}
\newcommand{\dm}{\diamondsuit}

\newcommand{\xto}{\xrightarrow}

\newcommand{\lra}{\leftrightarrow}

\newcommand{\RM}{\backslash}

\newcommand{\into}{\hookrightarrow}

\newcommand{\congto}{\xrightarrow{\sim}}
\renewcommand{\cong}{\simeq}

\newcommand{\bea}{\begin{eqnarray*} }
\newcommand{\eea}{\end{eqnarray*} }
\newcommand{\bee}{\begin{eqnarray} }
\newcommand{\eee}{\end{eqnarray} }
\newcommand{\be}{\begin{equation} }
\newcommand{\ee}{\end{equation} }
\newcommand{\bA}{\left(\begin{array}}
\newcommand{\eA}{\end{array}\right)}
\newcommand{\bma}{\begin{bmatrix}}
\newcommand{\ema}{\end{bmatrix}}
\newcommand{\bcd}{\begin{tikzcd}}
\newcommand{\ecd}{\end{tikzcd}}
\newcommand{\bcs}{\begin{cases}}
\newcommand{\ecs}{\end{cases}}

\newtheorem*{maintheo}{Main Theorem}
\newtheorem{theorem}{Theorem}[section]
\newtheorem{lemma} [theorem]{Lemma}
\newtheorem{proposition}[theorem] {Proposition}
\newtheorem{corollary} [theorem]{Corollary}
\newtheorem{remark} [theorem]{Remark}
\newtheorem{example}[theorem]{Example}
\newtheorem{definition}[theorem]{Definition}

\newcommand{\bp}{\begin{proposition}}
\newcommand{\ep}{\end{proposition}}
\newcommand{\bt}{\begin{theorem}}
\newcommand{\et}{\end{theorem}}
\newcommand{\bpf}{\begin{proof}}
\newcommand{\epf}{\end{proof}}
\newcommand{\bl}{\begin{lemma}}
\newcommand{\el}{\end{lemma}}
\newcommand{\bc}{\begin{corollary}}
\newcommand{\ec}{\end{corollary}}
\newcommand{\bd}{\begin{definition}}
\newcommand{\ed}{\end{definition}}
\newcommand{\bex}{\begin{example}}
\newcommand{\eex}{\end{example}}
\newcommand{\brem}{\begin{remark}}
\newcommand{\erem}{\end{remark}}
\newcommand{\bnum}{\begin{enumerate}}
\newcommand{\enum}{\end{enumerate}}

 \newcommand{\C}{\mathbb{C}}
  \newcommand{\D}{\mathbb{D}}

\newcommand{\R}{\mathbb{R}}
\newcommand{\Z}{\mathbb{Z}}

\newcommand{\ccal}{\mathcal{C}}

\newcommand{\hcal}{\mathcal{H}}

\newcommand{\lcal}{\mathcal{L}}

\newcommand{\scal}{\mathcal{S}}

\DeclareMathOperator{\Hom}{Hom}

\newcommand{\uhom}{\mathcal{H}om}

\renewcommand{\Im}{\text{Im}}

\DeclareMathOperator{\Mod}{Mod}

\newcommand{\orn}{\mathfrak{or}}

\DeclareMathOperator{\Perf}{Perf}

\newcommand{\tforall}{\text{ for all } \;}

\newcommand{\In}{\subset}

\newif\ifpaper
\papertrue

\usepackage{tikz, tikz-cd,pgfplots}
\usetikzlibrary{patterns, matrix,shapes,arrows,calc,topaths,intersections,
positioning,decorations.pathreplacing,decorations.pathmorphing,fit}
\newcommand*{\hrlen}{4}
\newcommand*{\hramp}{3}

\tikzset{
asdstyle/.style={blue,thick},
righthairs/.style={postaction={decorate,draw,decoration={border,amplitude=\hramp,segment length=\hrlen,angle=-90,pre=moveto,pre length=\hrlen/2}}},
lefthairs/.style={postaction={decorate,draw,decoration={border,amplitude=\hramp,segment length=\hrlen,angle=90,pre=moveto,pre length=\hrlen/2}}},
righthairsnogap/.style={postaction={decorate,draw,decoration={border,amplitude=\hramp,segment length=\hrlen,angle=-90}}},
lefthairsnogap/.style={postaction={decorate,draw,decoration={border,amplitude=\hramp,segment length=\hrlen,angle=90}}},
graphstyle/.style={thick},
arrowstyle/.style={thick,decorate,decoration={snake,amplitude=1.7,segment length=10pt,post length=.5mm,pre length=0}},
genmapstyle/.style={thick,-stealth'},
arrhdstyle/.style={thick},
exceptarcstyle/.style={red, ultra thick},
dualquiverstyle/.style={thick,->}
patstyle/.style={pattern color = gray, pattern = north east lines, opacity=0.3 }
}

\tikzset{
commutative diagrams/.cd,
arrow style=tikz,
diagrams={>=latex}}

\tikzset{->-/.style={decoration={
  markings,
  mark=at position #1 with {\arrow{>}}},postaction={decorate}}}
\tikzset{-<-/.style={decoration={
  markings,
  mark=at position #1 with {\arrow{<}}},postaction={decorate}}}

\newcommand{\linf}{\Lambda^\infty}
\newcommand{\La}{\Lambda}
\newcommand{\sinf}{SS^\infty}
\newcommand{\tinf}{T^\infty}
\newcommand{\lsm}{\Lambda_{sm}}

\title{Sheaf Quantization of Legendrian Isotopy}
\author{Peng Zhou}
\thanks{This work is supported by an
IHES Simons Postdoctoral Fellowship as part of the Simons
Collaboration on HMS. }

\address{Institut des Hautes \'Etudes Scientifiques. Le Bois-Marie, 35 route de Chartres, 91440 Bures-sur-Yvette France}
\email{pengzhou@ihes.fr}
\date{\today}

\begin{document}

\maketitle

\begin{abstract}
Let $\{\Lambda^\infty_t\}$ be an isotopy of Legendrians (possibly singular) in a unit cosphere bundle $S^*M$. Let $\ccal_t = Sh(M, \Lambda^\infty_t)$ be the differential graded (dg) derived category of constructible sheaves on $M$ with singular support at infinity contained in $\Lambda^\infty_t$. We prove that if the isotopy of Legendrians  embeds into an isotopy of Weinstein hypersurfaces, then the categories $\ccal_t$ are invariant. 
\end{abstract}

Let $M$ be a smooth compact manifold of real dimension $m$, $T^*M$ the cotangent bundle, and
\[\dot T^*M = T^*M - T^*_M M , \quad T^\infty M: =  \dot T^*M / \R_{>0}\] 
be the punctured cotangent bundle and contact cosphere bundle at infinity. Let $\linf \subset T^\infty M$ be a (singular) Legendrian, by which we mean a Whitney stratifiable subspace whose top dimensional strata are smooth Legendrian, and
\[ \dot \La = \R_{>0} \cdot \linf, \quad \La = \dot \La \cup T^*_M M \]
be two associated conical Lagrangians. 
We denote by 
$Sh(M, \Lambda^\infty)$ the the dg derived category of constructible sheaves  on $M$ where objects are sheaves $F$ with singular support at infinity $SS^\infty (F) \subset \Lambda^\infty$, i.e. $SS(F) \subset \La$. \footnote{One can work with either 'large', or 'traditional',  or 'wrapped' constructible sheaves \cite{N4}. Here for simplicity, we work with traditional constructible sheaf.}

%

\bd
Let $I \subset \R$ be an open interval, $(C, \xi)$ be a contact manifold. An {\bf isotopy of Legendrian} in $C$ over $I$ is a Whitney stratifiable closed subset $\lcal_I \subset C \times I$, such that $\lcal_t:=\lcal_I \cap C \times \{t\}$ is a (singular) Legendrian for all $t \in I$. We also denote an isotopy as $\{\lcal_t\}_{t \in I}$ or simply $\{\lcal_t\}$. 
\ed
\brem 
If we choose a contact form $\alpha$ on $(C,\xi)$, we may form a new contact manifold $(C \times T^* I, \alpha + \tau d t)$, and lift $\lcal_I$ to a Legendrian in $C \times T^*I$. The two description of isotopy are equivalent. 
\erem

We are interested in the following question: \\
\noindent{\bf  Main Question:} {\em Given an isotopy of Legendrians  $\{\Lambda^\infty_t\}$ in $T^\infty M$, when is it 'non-characteristic' \cite{N3}, that is, the sheaf category $Sh(M, \Lambda^\infty_t)$ remains invariant? Or more concretely, if we deform $\Lambda^\infty$, can we deform the sheaf $F$ such that $SS^\infty(F)$ remains in $\Lambda^\infty$?}


Before we state our main result, we first review two important results in this direction. The first one is due to Guillermou-Kashiwara-Schapira, which quantizes isotopy of the entire contact manifold $T^\infty M$. 
\bt[\cite{GKS} Theorem 3.7, Proposition 3.12] \label{t:GKS}
Let $I$ be an open interval containing $0$, and $\varphi: I \times T^\infty M\to T^\infty M$ be a smooth map with $\varphi_t=\varphi(t,-)$. Assume $\varphi$ satisfies (1) $\varphi_0=id$, and (2) $\varphi_t$ are contactomorphisms for all $t \in I$. Then for each $t \in I$, we have equivalences of category
\[ \h \varphi_t: Sh(M) \congto Sh(M), \quad \z{ such that } \; SS^\infty(\h \varphi_t F) = \varphi_t (SS^\infty(F)). \]
\et

One immediately get the following corollary. 
\bc
If the isotopy of Legendrian $\{\Lambda^\infty_t\}_{t \in I} $ can be embedded into an isotopy $\{\varphi_t\}_{t \in I}:S^*M \to S^*M$ of the contact manifold, that is,  $\Lambda^\infty_t = \varphi_t(\Lambda^\infty_0)$. Then we have equivalence of categories
\[ \h \varphi_t: Sh(M, \linf_0) \congto Sh(M, \linf_t). \]
\ec
\brem
Any isotopy of smooth Legendrian can be extended to a contact isotopy of the ambient manifold. 
If the Legendrian is singular, and if the homeomorphism type of the Legendrian changes during the isotopy, then it cannot be extended to a contact isotopy.  

\ifpaper 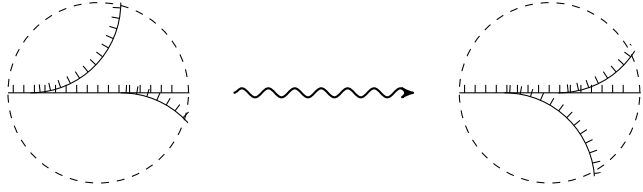
\begin{figure}[h]
\centering
%
%
\begin{tikzpicture}[scale=0.6]
\begin{scope}
   \clip (0,0) circle (2.1);  \draw [dashed] (0,0) circle (2); 
\draw [lefthairs] (-2,0) -- (2,0); 
\draw [lefthairs] ( 0.5,0) arc (90:0:2);
\draw [righthairs] (0.5,2) arc (0:-90:2);
\end{scope}


\begin{scope} [ shift = {(10, 0)} ]
 \clip (0,0) circle (2.1);  \draw [dashed] (0,0) circle (2); 
\draw [lefthairs] (-2,0) -- (2,0); 
\draw [lefthairs] (-1,0) arc (90:0:2);
\draw [righthairs] (2.2,2) arc (0:-90:2);
\end{scope}

\draw [arrowstyle, -stealth'] (3,0) -- (7,0); 
\end{tikzpicture}
\caption{An example of Legendrian isotopy (shown as front projection from $S^*\R^2 \to \R^2$, with `hairs' indicating co-direction) which cannot be embedded in a contact isotopy. }
\end{figure} \fi

\erem

The second result is due to Nadler \cite{N3}, where he proves that any Legendrian singularity admits a non-characteristic deformation to an arboreal singularity (introduced in \cite{N2}). In  \cite{N3}, Nadler proposed the following geometric condition on  Legendrian isotopies. 
\bd  [Displaceable Legendrian]
Let $(\tinf M, \xi = \ker \alpha, R_\alpha)$ be the cosphere bundle with a choice of Reeb vector field $R_\alpha$, and let $R_\alpha^t: \tinf M \to \tinf M$ be the Reeb flow for time $t$.\footnote{Throughout the paper, we will use the notation $X^t$ for the flow generated by a vector field $X$ for time $t$.}   Let $\epsilon>0$.  A Legendrian $\linf \subset \tinf M$ is {\bf displaceable} for $R_\alpha$ if there exists a constant $\epsilon>0$, such that 
\be \label{e:sep} \Lambda^\infty  \cap R_\alpha^s(\Lambda^\infty) = \emptyset, \quad \forall 0 < |s| < \epsilon.  \ee

We say a family of Legendrian $\{\linf_t\}$ is {\bf uniformly displaceable} for $R_\alpha$, if each $\linf_t$ is displaceable for the same constant $\epsilon$. 
\ed

 \ifpaper \begin{figure}[H]
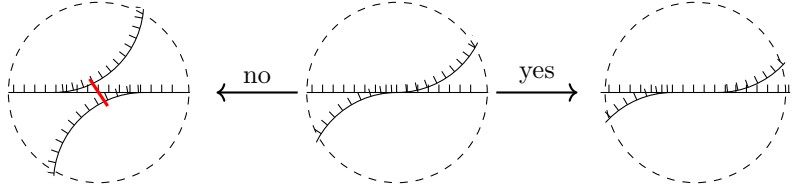

\tikz[scale=0.6]{
\clip (0,0) circle (2.1);  \draw [dashed] (0,0) circle (2); 
\draw [lefthairs] (-2,0) -- (2,0); 
\draw [righthairs] (1,0) arc (90:180:2);
\draw [righthairs] (1,2) arc (0:-90:2);
\draw [red, very thick] (0.2, -0.3) -- (-0.2, 0.3); 
}
\,
\tikz[scale=0.6]{
\begin{scope}
\clip (0,0) circle (2.1);  \draw [dashed] (0,0) circle (2); 
\draw [lefthairs] (-2,0) -- (2,0); 
\draw [righthairs] (0,0) arc (90:180:2);
\draw [righthairs] (2,2) arc (0:-90:2);
\end{scope}
\draw [thick,->] (2.2,0) -- node[above,pos=0.5] {yes}  (4,0);
\draw [thick,->] (-2.2,0) -- node[above,pos=0.5] {no} (-4,0);
}
\,
\tikz[scale=0.6]{
\clip (0,0) circle (2.1);  \draw [dashed] (0,0) circle (2); 
\draw [lefthairs] (-2.5,0) -- (2.5,0); 
\draw [righthairs] (-0.5,0) arc (90:180:2);
\draw [righthairs] (2.5,2) arc (0:-90:2);
}
\caption{The deformation  to the right is uniformly displaceable, and the one to the left is not, due to the appearance of new short Reeb chord (marked in red). (c.f. \cite{N3}, Example 1.5)}
\end{figure} \fi

It turns out just having the uniform displaceablity for Legendrian is not enough, one need to impose some control on the topology of the Legendrians as well. 
\bex
Let $C = J^1 \R, \alpha = dz - y dx, R_\alpha = \pa_z$. For $t \in [0,1)$, define a family of Legendrians
\[ \lcal_t = \{(x, 0,0): x \in \R\} \cup \{(x, x^2+t,x^3/3 + tx): x \in \R\} \]
there are no Reeb chords ending on $\lcal_t$ even at $t=0$. However the sheaf category $\ccal_t$ associated to $\lcal_t$ (after identifying $J^1\R$ with $S^*\R^2$ with matching contact forms) jumps as $t \to 0^+$, since the topology of the Legendrian changed. 

\ifpaper 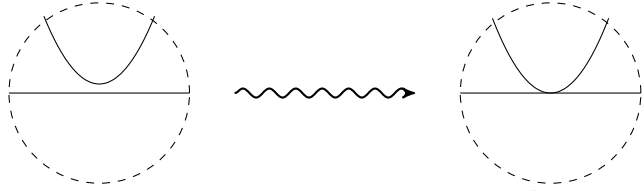
\begin{figure}[H]
\begin{tikzpicture}[scale=0.6]
\begin{scope}
   \clip (0,0) circle (2.1);  \draw [dashed] (0,0) circle (2); 
\draw  (-2,0) -- (2,0); 
\draw  plot[smooth,domain=-2:2] (\x, {0.2+\x*\x});
\end{scope}


\begin{scope} [ shift = {(10, 0)} ]
\clip (0,0) circle (2.1);  \draw [dashed] (0,0) circle (2); 
\draw  (-2,0) -- (2,0); 
\draw  plot[smooth,domain=-2:2] (\x, {\x*\x});
\end{scope}

\draw [arrowstyle, -stealth'] (3,0) -- (7,0); 
\end{tikzpicture}
\caption{The deformation of Legendrian $\lcal$ in $J^1\R \cong T^*\R \times \R$, draw as the Lagrangian projection to $T^*\R$. Note there is no Reeb chord ending on $\lcal$ throughout the deformation. \label{noreeb}}
\end{figure} \fi
\eex

\ss{Definitions and Result}


To state our main result, we need some definitions. Recall that a hypersurface in a contact manifold is {\em convex} \cite{Gi} if it admits a transverse contact vector field. We want to consider nested tubular neighborhoods of Legendrians with convex boundaries. 

\bd\label{d:tube}
Let $\lcal$ be a singular Legendrian in $(C, \xi)$. 
A  {\bf convex tubular neighborhood}  for $\lcal$ is the following data $(U, \rho, X)$, \\
(1) $U=:U(\lcal)$ is an open tubular neighborhood of $\lcal$, \\
(2) $\rho:U \to [0,1)$ is a $C^1$-function,\\
(3) $X$ is a smooth contact vector field on $U$. \\
Such that if we define 
\[ U_r(\lcal) = \{x \in U: \rho(x)<r\}, \quad \forall 0<r<1 \] 
then,  $\cap_r U_r(\lcal) = \lcal$ ; $\pa U_r(\lcal)$ are $C^1$-smooth and $C^1$-diffeomorphic; $X$ is transverse to all $\pa U_r(\lcal)$; and $d\rho(X) > c \rho$ for some constant $c>0$. 

Let $\{\lcal_t\}$ be an isotopy of Legendrian in $(C, \xi)$. 
An {\bf isotopy of convex tubular neighborhoods} $\{(U,\rho,X)_t\}$ of $\{\lcal_t\}$ is a one-parameter family of such data with uniform bound on constant. 
\ed

The notion of Weinstein hypersurface is introduced in \cite{Av}. We recall its definition following \cite[Section 2]{Eli}. 
\bd \label{d:hyp} 
\begin{enumerate}
\item  A codimension-$1$ submanifold $\hcal$ in a contact manifold $(C, \xi)$ with boundary $\pa \hcal$ is called {\bf Weinstein hypersurface} if there exists a contact form $\alpha$ such that $(\hcal, \lambda:=\alpha|_\hcal)$ is compatible with a Weinstein structure on $\hcal$, i.e. $d\lambda$ is symplectic and the Liouville vector field $X$ on $\hcal$ dual to the Liouville form $\lambda$ is outward transverse to $\pa \hcal$ and admits a Lyapunov function $\phi: \hcal \to \R$.  

We denote a Weinstein hypersurface (including a specification of compatible Weinstein structure) by $(\hcal, \lambda, X, \phi)$. 
\item If $\lcal$ is the skeleton of $(\hcal, \lambda, X, \phi)$, we say $\hcal$ is a {\bf Weinstein hypersurface thickening} of $\lcal$. 
\item An isotopy of Weinstein hypersurface is a smooth family of $\{(\hcal, \lambda, X, \phi)_t\}$ where the choice of contact 1-form $\alpha_t$ has smooth and bounded variation with $t$. 
\end{enumerate}
\ed

We can show that if the Legendrian admits a Weinstein hypersurface thickening and is uniformly displaceable, then it admits a canonical tube thickening (Proposition \ref{p:hyp2tube}). 

\begin{maintheo} [Theorem \ref{t:ext}]
Let $\{\linf_t\}_{t \in I}$  be an isotopy of Legendrian in $\tinf M$, such that $\linf_t$ is constant for $t$ outside a closed interval $[a,b] \subset I$. If $\{\linf_t\}$ is uniformly displaceable for some Reeb vector field $R_\alpha$ on $\tinf M$,  and
\begin{enumerate}
\item there exists an isotopy of convex tubular neighborhoods $\{(U(\linf_t), \rho_t, X_t)\}$ of $\{\linf_t\}$
\item or, there exists an isotopy of Weinstein hypersurface thickening $\{(H, \lambda, X, \phi)_t\}$ of $\{\linf_t\}$,
\end{enumerate}
then the sheaf categories $Sh(M, \linf_t)$ remain invariant.
\end{maintheo}

\brem
The notion of convex tubular neighborhood is related to the `frozen boundary' for a Liouville sector \cite{GPS}, and the Weinstein hypersurface is related to `stop' in partially wrapped Fukaya category \cite{Syl}. See \cite[Section 2,3]{Eli} for related work on Weinstein hypersurface and Weinstein pair. 
\erem

\ss{Idea of the Proof}
We first give an heuristic derivation for why we might expect such a theorem, though we do not follow this approach literally. See the previous section for notations $X_s, U_r(\lcal_s), \cdots $. The main idea is to use the retracting contact flow $-X_s$ toward $\linf_s$, properly cut-off outside $U_{\delta}(\lcal_s)$ for some $1/2<\delta<1$,  to deform and squeeze a nearby Legendrian skeleton $\linf_t \subset U_{\delta}(\linf_s)$ into $\linf_s$ in the limit. We consider the sheaf quantization of the retracting flow for time $T>0$, 
\[ X_s^{-T}: \tinf M \congto \tinf M \quad \leadsto \quad \h X_s^{-T}: Sh(M) \congto Sh(M). \]
Then we define the projection functors as the limit of the flow
\be \Pi_s: Sh(M, U_\delta(\linf_s)) \to Sh(M, \linf_s), \quad
 \Pi_s(F) := \lim_{T \to \infty} \h X_s^{-T} (F_t), \label{Pis}\ee
 where $Sh(M, U_\delta(\linf_s))$ means constructible sheaves with $SS^\infty(F) \subset U_\delta(\linf_s)$. The limit is not inductive, or projective limit, and is defined (in the style of a nearby cycle functor) in Section \ref{limcon} . Then, one only need to show that for any $t,s$ closed enough (contained in each other's tubular neighborhoods $U_\delta$), we have a pair of inverse functors 
 \[ \Pi_s|_{Sh(M, \linf_t)}: Sh(M, \linf_t) \xleftrightarrow{\sim} Sh(M, \linf_s) : \Pi_t|_{Sh(M, \linf_s)}. \]
To see they are inverses, we consider a constructible sheaf $F_t \in Sh(M, \linf_t)$ as functors $\Hom(-, F_t)$, and test on 'probe' sheaves $P$ such that 
\be \sinf(P) \cap [U_\delta(\linf_t) \cup U_\delta(\linf_s)] = \emptyset \label{disjP}\ee 
then $\Hom(P, \h X_t^{-T_1} \h X_s^{-T_2} F_t)$ is independent of $T_1 ,T_2$. One way to construct such probe uses {\em wrapped constructible sheaves} (see \cite{N4} for definition), we have \cite[Theorem 1.6] {N4}
\[ Sh(M,  \linf_t) \congto \z{Fun}^{ex}(Sh^w(M, \linf_t)^{op}, \z{Perf}_\C). \]
To achieve \eqref{disjP}, we use small {\em negative} Reeb flow to displace $P$ without changing the homs (Proposition \ref{p:nc2}). We get, for all $P \in Sh^w(M, \linf_t), F_t \in Sh(M, \linf_t)$
\bea \Hom(P, F_t) \cong \Hom(\h R^{-\epsilon} P, F_t) \cong \Hom(\h R^{-\epsilon} P, \Pi_s(F_t)) \\
\cong \Hom(\h R^{-\epsilon} P, \Pi_t \Pi_s(F_t)) \cong \Hom(P, \Pi_t \Pi_s(F_t))) \eea
Hence $\Pi_t, \Pi_s$ are inverses. 
\vskip 5mm

Our actual approach is as following: let $\linf_I \subset \tinf(M \times I)$ be an isotopy of Legendrians and let $F_t \in Sh(M, \linf_t)$. We will extend $F_t$ to a sheaf $F_I \subset Sh(M \times I, \linf_I)$ such that $F_I|_t=F_t$. 

One first show that such extension is unique (if exists), this is equivalent to show that restriction functor $F_I \mapsto F_t$ is fully-faithful, i.e. (Proposition \ref{p:uniq})
\[ \Hom(F_I, G_I) \congto \Hom(F_t, G_t), \forall F_I, G_I \in Sh(M \times I, \linf_I). \]
One need to show that $\uhom(F_I, G_I)(M \times (a,b))$ is independent of the size of the interval, hence one can interpolate from $(a,b) = I$ to infinitesimal small neighborhood around $t$. The key technical point is to use the uniform displaceability condition to perturb $G_I$ slice-wise by positive Reeb flow for time $s$, $G_I \to K_s^!G_I$, to separate $\sinf(F_I)$ and $\sinf(K_s^!G_I)$. 

One then show that such extension exists locally, i.e., given $F_t$, we may find a small neighborhood $(t-\delta, t+\delta)$, where $\delta$ is uniform, to extend $F_t$ on $M \times \{t\}$ to $M \times (t-\delta, t+\delta)$. This is done using limit of the retracting flow, as done in defining $\Pi_s$ in \eqref{Pis}. For general contact flow, there is no way to take limit. Here we can take limit since the singular support of the sheaf $\sinf(F_t)$ converges under $-X_s$ to the sink of the flow $\linf_s$. We thus get the limiting sheaf with desired bound on singular support. 

Finally, we use uniqueness of extension to patch together local extensions, and get the global extension result. (c.f. Lemma 1.13 in \cite{GKS}). 

\brem
We thank V. Shende for informing us the up-coming work of Nadler-Shende about quantization of exact symplectic category, which include a result on invariance of microlocal sheaf category $Sh(W)$ for Weinstein manifold $(W, \lambda)$ \cite{Sh} under Weinstein homotopy.
\erem

\ss{Acknowledgements}
I would like to thank my advisor Eric Zaslow for suggesting the idea of 'invariance of hom under Reeb perturbation'. I also thank P. Schapira for many warm discussions, and for suggesting finding a sheaf-theoretic proof for Proposition \ref{p:nc2}.  I thank D. Nadler for encouragements and comment on an early draft of this paper using almost retraction. I also thank S. Guillermou for explaining many points in the \cite{GKS} paper,  and V. Shende for many useful discussions. 

\section{Convex Tubular Neighborhoods and Weinstein Hypersurfaces \label{s:thicken}}
We give basic definition and construction for Weinstein hypersurface and convex tubular neighborhood. We will work with general contact manifold $(C, \xi)$ instead of $S^*M$ so that the results may generalize to other Weinstein domain.  

\ss{Basic of Contact Geometry}\label{contact}
We recall the definition of co-oriented contact manifold as follow. Let $C$ be a $2n+1$ dimensional manifold, $\xi \subset TC$ be a rank $2n$ sub-bundle, such that there exists a one-form (contact one-form) $\alpha$ (up to multiplication of non-negative function) satisfying $\xi = \ker \alpha$ and $\alpha \wedge (d\alpha)^n \neq 0$. If we fix such a $\alpha$, we have a Reeb vector field $R_\alpha$ given by 
\[ \iota_{R_\alpha} \alpha = 1, \quad \iota_{R_\alpha} d\alpha = 0. \]
We note that different choices of $\alpha$ will lead to different choices of $R_\alpha$. 

A contact vector field $X$ is one that perserves $\xi$. 
\bd
Given a smooth function $H: C \to \R$, the {\bf contact Hamiltonian vector field} $X_H$ is uniquely determined by
\be \bcs \la X_H, \alpha \ra = H  \\ 
\iota_{X_H} d \alpha =  \la H, R \ra \alpha - d H \ecs  \ee
\ed
Reeb vector field is a speical case of $X_H$ for $H=1$. 

\bp[\cite{Ge} Theorem 2.3.1]\label{p:XH}
With a fixed choice of contact form $\alpha$ there is a one-to-one correspondence between contact vector field $X$ and smooth functions $H: C \to \R$. The correspondence is given by
\[ X \mapsto H = \la \alpha, X \ra, \quad H \mapsto X_H. \]
\ep

Unlike symplectic Hamiltonian vector field, $X_H$ does not preserve level set of $H$. 
\bl \label{l:XHH}
\[ \la X_H, dH \ra = H \la R, dH \ra \]
In particular, $X_H$ preserves the zero set of $H$. 
\el
\bpf
Since $X_H = H R + X_H^\|$, where $X_H^\| \in \ker(\alpha)$, we have 
\[ \la X_H - H R, d H \ra = \la X_H - H R, d H - \alpha \ra =  \la X_H - H R, -\iota_{X_H}(d\alpha) \ra \]\[= d\alpha(X_H - H R, X_H) = 0 \]
where we have used $R \in \ker(d\alpha)$. 
\epf

\bex
Let $M$ be a smooth manifold, and $T^*M$ the cotangent bundle with canonical Liouville one-form $\lambda$ and symplectic two-form $\omega=d\lambda$. If we put local Darboux coordinate $(q,p)=(q_1,\cdots, q_m; p_1, \cdots, p_m)$ on $T^*M$ where $m=\dim_\R M$, then $\lambda = \sum_{i=1}^m p_i dq_i$ and $\omega = \sum_i dp_i \wedge dq_i$, and we will suppress the indices and summation to write $\lambda = pdq, \omega = dpdq$.  Also define $\dot{T}^*M = T^*M \RM T^*_M M$, $T^\infty M = \dot T^*M / \R_{>0}$. The Liouville vector field for $\lambda$ is defined by defined by $\iota_{V_\lambda} \omega = \lambda$, and here it is given by $V_\lambda = p \d_p$.  On $T(\dot T^*M)$, the symplectic orthogonal to the Liouville vector field defines a distribution 
\[ \wt \xi = \{(q,p; v_q, v_p) \in T(\dot T^*M): \omega((v_q, v_p), V_\lambda) = 0\}, \] which project to a canonical contact distribution $\xi$ on $T^\infty M$. Let $g$ be any Riemaninan metric on $M$, then $T^*M$ has induced norm.  Let $S^*M = \{ (q, p) \in T^*M \mid |p|=1\}$ be the unit cosphere bundle with contact form $\alpha = \lambda|_{S^*M}$, then the contact distribution can also be written as $\xi = \ker (\alpha)$. 
\eex

Define the symplectization of $(C, \xi = \ker \alpha)$ by 
\[ S:= C \times \R_u, \quad \lambda = e^u \alpha, \quad \omega_S = d (e^u \alpha). \]
We have projection along $\R_u$, and inclusion of zero section as: 
\[ \pi_S: S \to C, \quad \iota_C: C \cong C \times \{0\} \into S. \]
Different choice of $\alpha$ gives the same $S$ up to fiber preserving symplectomorphism, that identifies the 'zero-section' $\Im (\iota_C)$. 

A Hamiltonian function $H: C \to \R$ can be extended to a homoegeneous degree one function $\wt H: S \to \R$ by $\wt H =e^u H$. Then the symplectic Hamiltonian vector field  $\xi_{\wt H}$, given by $\omega_S(-, \xi_{\wt H}) = d \wt H(-)$, preserves the fiber of $\pi_S$ and descend to $X_H$. 

\subsection{Weinstein Hypersurface}
Let $(\hcal, \lambda, X,\phi)$ be a Weinstein hypersurface in $(C,\xi=\ker(\alpha))$, in particular $\lambda = \alpha|_\hcal$ (Definition \ref{d:hyp}).  
\footnote{ 
If we only fix $\hcal$ but allowing $\lambda$ and $\alpha$ to vary, then we may change $\alpha$ to $e^f \alpha$ for some smooth function $f$ as long as $k:=1+ \la df, X \ra>0$ on $\hcal$, in this case the Liouville field changes to $\frac{1}{k}X$, and is gradient-like for the same $\phi$. Moreover, the skeleton for $\hcal$ remains the same. See \cite[Section 2] {Eli} and \cite[Lemma 12.1]{CE}. 
}
For small enough $\epsilon>0$, we may thicken $\hcal$ to $U_\epsilon(\hcal)$ by Reeb flow for time $|t|<\epsilon$.
We will take contact Hamiltonian function $H = t$ on $U_\epsilon(\hcal)$, where $t$ is the time coordinate (not to be confused with the isotopy parameter later). Then 
\[ \la R, dH \ra = \la \pa t, dt \ra = 1, \quad \hcal = \{H=0\}. \]
 
\bp\label{p:XH}
Under the identification $U_\epsilon(\hcal) \simeq \hcal \times (-\epsilon, \epsilon)$, the contact form $\alpha$ and Reeb vector field can be written as 
\[ \alpha = \lambda + d t,\quad R = \pa_t \]
The contact vector field $X_H$ can be written as 
\[ X_H = X + t \pa_t, \]
where $\lambda$ is the Liouville form on $\hcal$, $X$ the Liouville vector field along $\hcal$.  
\ep
\bpf
We call $\hcal$ the horizontal direction and $(-\epsilon, +\epsilon)$ the vertical direction. Since the Reeb flow is translation the $t$ coordinate, and the Reeb flow preserves $\alpha$,  we have $\alpha = \lambda$ on the  horizontal direction. Since $\iota_R \alpha=1$, we have $\alpha = dt$ along the vertical direction. Thus $\alpha = \lambda + d t$. We note that $\iota_X \lambda = \iota_X(\iota_X(d\lambda)) = 0$, and $\iota_X(dt)=0$, hence $X \in \ker \alpha$. Thus we may easily check the given formula $X_H$ satisfy the definition. 
\epf

\bc\label{mfd2hyp}
Let $(\hcal, \lambda, X, \phi) \into (C ,\xi = \ker \alpha)$ be a Weinstein hypersurface. If $\{(\hcal, \lambda_\delta, X_\delta, \phi_\delta)\}_{|\delta|<c}$ is an isotopy of Weinstein domain, where 
\[ \lambda_\delta = \lambda + \delta d f \] 
for some smooth and uniformly bounded $f: \hcal \to \R$. Then there exists an isotopy of Weinstein hypersurface realizing the given isotopy of Weinstein domain for $|\delta| < c'$ where $c' < c$. 
\ec
\bpf
We work in the neighborhood $U_\epsilon(\hcal) \cong \hcal \times (-\epsilon, \epsilon)$ with coordinate $(x,t)$. Since $\alpha = \lambda + dt$ in $U_\epsilon(\hcal)$, we may define a family of hypersurface as graph of $\delta f$
\[\hcal_\delta = \{(x,t) \in \hcal \times (-\epsilon, \epsilon) \mid t= \delta f(x) \}, \forall |\delta| < c'\]
where 
\[  c':= \max\{\delta \mid \delta < c, \sup_{x \in \hcal} \delta f(x) < \epsilon \}. \]
We have canonical identification $\pi_\delta: \hcal_\delta \to \hcal$ by project away the $t$ coordinate, and $\lambda_{\hcal_\delta}: = \alpha|_{\hcal} = \lambda + \delta df$.  
\epf

\ss{Construction of Convex Tubular Neighborhood}

\bp \label{p:hyprho}
Let $(\hcal, \lambda, X, \phi)$ be a Weinstein domain, such that $\pa \hcal = \phi^{-1}(c)$. There exists a unique $C^1$-function $\psi: \hcal \to [0,1]$, such that $\psi$ is smooth on $\hcal \RM (\pa \hcal \cup \lcal)$, $\psi|_{\pa \hcal} = 1$, $\psi|_{\lcal}=0$, and $\la X, d \psi \ra = 2\psi$ on $\hcal \RM \lcal$.  
\ep
\bpf
We have a diffeomorphism generated by the downward Liouville flow $-X$
\[ \Psi: \pa \hcal \times \R_{\geq 0} \congto \hcal \RM \lcal, \quad (x, t) \mapsto X^{-t}(x). \]
Thus we may define $\psi$ on $\hcal \RM \lcal$ by 
\[  \psi|_{\hcal\RM \lcal} (\Psi^{-1}(x,t)) = e^{-2t}.\]
Thus  $-\pa_t e^{-2t} = \la d\psi, X \ra =2 \psi$.  $\psi$ has a $C^1$ extension by zero to $\lcal$, since $d \psi |_\lcal=0$. 
\epf

Let $U=U_\epsilon(\hcal) \simeq \hcal \times (-\epsilon, \epsilon)$ with $x$ coordinate on $\hcal$ and $t$ coordinate on $(-\epsilon, \epsilon)$. We define
\be \label{e:Urho}   \rho(x,u) = \psi(x) + t^2: U \to \R. \ee
Then 
\bp\label{p:Urho}
$  \rho$ is a $C^1$-function on $U$, vanishing only on $\lcal$;  and $\la d   \rho, X \ra =2   \rho > 0$ on $U \RM \lcal$. 
\ep
\bpf
The regularity and vanishing statement is clear. Using Proposition \ref{p:XH} and \ref{p:hyprho}, we have 
\[ \la d(\psi + t^2), X_H \ra = \la d\psi + 2t dt, X + t\pa_t \ra = 2\psi + 2t^2 = 2  \rho \]
away from $\lcal$. 
\epf

\bp
If $(\hcal, \lambda, V, \phi)$ is a Weinstein hypersurface thickening of Legendrian $\lcal$, then there exist a convex tubular neighborhood thickening $(U, \rho, X)$ of $\lcal$. 
\ep
\bpf
Let $\alpha$ be the family of contact 1-form, such that $\lambda = \alpha|_{\hcal}$. Let $1 \gg \epsilon>0$ be small enough, such that 
\[ R_{\alpha}^{\tau} (\hcal) \cap \hcal = \emptyset, \quad \forall 0< |\tau| < 2 \epsilon. \]
Then we define $U^0 = U_{\epsilon}(\hcal)$ and $H$ using 1-form $\alpha$ and the associated Reeb flow.   Let $\rho^0$ be the $C^1$-function $\rho$ as constructed in \eqref{e:Urho}. We see $U^0$ and $\rho^0$ depends on $(\hcal, \lambda, X, \phi)$ and $\alpha$ canonically, hence $U^0$ has piecewise smooth boundary 
\[ \pa U^0 = \pa \hcal \times (-\epsilon, \epsilon) \cup \hcal \times \{-\epsilon, \epsilon\} \]
and $\rho^0$ is a globally $C^1$-function defined on $U^0$ and it is smooth away from $ \lcal \times (-\epsilon, \epsilon).$ The vector field $X_{H}$ is smooth in $U$ and is smoothly varying in $t$. 

Since $\{ \rho^0(x) < \epsilon^2 \}$ is contained in $U^0$, we can trim $U^0$ and rescale $\rho^0$ by
\[ U :=\{ \rho^0(x) < \epsilon^2 \}, \quad \rho:= \rho^0/\epsilon^2. \]
Let $X := X_{H}$, we still have 
\[ \la X, d\rho \ra =  \epsilon^{-2} \la X, d\rho^0 \ra = 2 \epsilon^{-2} \rho^0 = 2 \rho. \]
Then the data $(U, \rho, X)$ forms a convex tubular neighborhood thickening $\lcal$. 
\epf

\bp \label{p:hyp2tube}
Assume $\{\lcal_t\}_{t \in I}$ is an isotopy of Legendrian, uniformly displaceable for some Reeb vector field, and  can be thickened to a Weinstein hypersurface isotopy $\{(\hcal, \lambda, V, \phi)_t\}$. Then there exist a convex tubular neighborhood thickening $\{(U, \rho, X)_t\}. $
\ep
\bpf
Since all the parameters have smooth and bounded dependence on $t$ (needed if $I$ is not compact), hence the proof of Proposition \ref{p:hyp2tube} goes through verbatim. 
\epf

\section{Non-Characteristic Isotopy of Sheaves}

\ss{Constructible Sheaves}

We give a quick working definition for constructible sheaf used here, and point to \cite{KS, S} for proper treatment. A constructible sheaf $F$ on $M$ is a sheaf valued in chain complex of $\C$-vector spaces, such that its cohomology is locally constant with finite rank with respect to some Whitney stratification $\scal = \{\scal_\alpha\}_{\alpha \in A}$ on $M$, where $\scal_\alpha$ are disjoint locally closed smooth submanifolds with nice adjacency condition and  $M= \sqcup_{\alpha \in A} \scal_\alpha$. The singular support $SS(F)$ of $F$ is a closed conical Lagrangian in $T^*M$, contained in $\cup_{\alpha \in A} T_{\scal_\alpha}^* M$, such that $SS(F) \cap T^*_M M$ equals the support of $F$,  and $(p,q) \in SS(F) \RM T^*_M M$ if there exists a locally defined function $f$ with $f(q)=0, df(q)=p$, such that the restriction map 
$ F(B_\epsilon(q) \cap \{f<\delta\}) \to F(B_\epsilon(q) \cap \{f<-\delta\}) $
fails to be a quasi-isomorphism for $0 < \delta \ll \epsilon \ll 1$. We denote by $SS^\infty(F) = SS(F) \cap S^*M$ the singular support of $F$ at infinity.

If $\Lambda \subset T^*M$ is a conical Lagrangian containing zero section (as always), we write $Sh(M, \Lambda^\infty)$ for the dg derived category of constructible sheaves \cite{N1} with object $F$ satisfying $SS^\infty(F) \subset \Lambda^\infty$. 

\bex
For example,  on $\R$, if $\C_{[0,1]}$ (resp. $\C_{(0,1)}$) denote constant sheaf with stalk $\C$ on $[0,1]$ (resp. on $(0,1)$) and zero stalk elsewhere, then their singular supports in $T^*\R$ are
\[  SS(\C_{[0,1]}) = 
\tikz[baseline=-1] {
\draw [gray, opacity=0.3] (-0.5, -0.8) grid (1.5, 0.8); 
\draw [very thick, red](0,0.8) -- (0,0) -- (1,0) -- (1,-0.8); 
},
\quad
SS(\C_{(0,1)}) = 
\tikz[baseline=-1] {
\draw [gray, opacity=0.3] (-0.5, -0.8) grid (1.5, 0.8); 
\draw [very thick, red] (0,-0.8) -- (0,0) -- (1,0) -- (1,0.8); 
}. 
\]
\eex

\bex
Let $j: U=B(0,1) \into \R^2$ be the inclusion of an open unit ball in $\R^2$. Then $j_* \C_U$ is supported on the closed set $\wb U$, with singular support at infinity as 
\[ SS^\infty(j_* \C_U) = \{ (x, \eta) \in S^*\R^2 \mid x \in \pa U, \eta = -d |x| \} = \quad \tikz[baseline=-3pt] {\draw [lefthairs] (0,0) circle (0.5);} \]
And $j_! \C_U$ is supported on the open set $U$, with singular support at infinity as
\[ SS^\infty(j_! \C_U) = \{ (x, \eta) \in S^*\R^2 \mid x \in \pa U, \eta = d |x| \} = \quad \tikz[baseline=-3pt] {\draw [righthairs] (0,0) circle (.5);}\]
Here  the Legendrians are represented by co-oriented hypersurfaces in $\R^2$ with hairs indicating the co-orientation. 
\eex

\ss{Operation on Constructible Sheaves}
Let $X, Y$ be manifolds. We use $f^*, f_*, f^!, f_!, \uhom, \ot$ to mean the corresponding dg derived functors: 
\bea
  - \ot F: & Sh(X) \lra Sh(X) & : \uhom(F, -)\\
  f^*:& Sh(X) \lra Sh(Y) &: f_* \\
  f_!: &Sh(Y) \lra Sh(X)&: f^! 
\eea
where $f: Y \to X$ is a map of real analytic manifolds. 

The Verdier duality $\D: Sh(X) \to Sh(X)$ is an anti-involution. It interchanges shriek with star
\[ \D\D = id, \quad f_! = \D f_* \D, \quad f^! = \D f^* \D. \]
The shrieks and stars are directly related in two cases: when $f$ is proper $f_! = f_*$; when $f$ is a smooth morphism of relative dimension $d_f$, $f^!(-) \cong f^*(-)\otimes  \omega_{Y/X} \cong f^*(-)\otimes \orn_{Y/X} [d_f]$, where $\orn_{Y/X}$ is the orientation sheaf of the fiber.

Given an open subset $U$ of $X$ and its closed complement $Z$,  
\[  \text{open inclusion:} \quad U \xhookrightarrow{j} X \xhookleftarrow{i} Z, \quad \text{closed inclusion},  \]
we have $j^* = j^!$ and $i_*=i_!$. Furthermore, there are exact triangles
\[ i_! i^! \to id \to j_* j^* \xto{[1]}, \quad  j_! j^! \to id \to i_* i^* \xto{[1]}.  \]
These are sheaf-theoretic incarnations of excisions: applied to the constant sheaf on $X$ and taking global sections, we get
\[ H^*(Z, i^! \C) \to H^*(X, \C) \to H^*(U, \C) \xto{[1]}, \quad H^*_c(U, \C) \to H^*_c(X, \C) \to H^*_c(Z, \C) \xto{[1]}. \]

If $Y$ is a locally closed $\ccal$-submanifold of $X$, we use $j_Y: Y \into X$ to denote the inclusion. Let $\C_Y \in Sh(Y)$ denote the constant sheaf on $Y$, and $\omega_Y = \D \C_Y$ be the Verdier dualizing complex of $Y$, then $\omega_Y$ is the canonically isomorphic to the shifted orientation sheaf $ \orn_Y[\dim Y]$  on $Y$. The standard sheaf on $Y$ is $j_{Y*} \C_Y$, and the costandard sheaf on $Y$ is $j_{Y!} \omega_Y$.

Let $X_i$, $i=1,2$, be spaces, and $K \in Sh(X_1 \times X_2)$. We define the following pair of adjoint functors
\be K_!: Sh(X_1) \lra Sh(X_2): K^! \ee
\be K_!: F \mapsto {\pi_2}_! (K \ot \pi_1^* F), \qquad  K^!: G \mapsto {\pi_1}_*(\uhom(K, \pi_2^! G)) \label{K:shriek}\ee
In \cite{KS},$K_! = \Phi_K$ and $K^!=\Psi_K$ and with $X_1, X_2$ switched. The notation here is suggestive for them to be adjoint functors. 

\ss{Isotopy of Legendrian and Sheaves}
Let $I =(0,1) \subset \R$. For any $t \in I$, let 
\[ j_t: M_t := M \times \{t\} \into M_I := M \times I \] 
be the inclusion of $t$-slice $M_t$ into the total space $M_I$, and let $\pi_I: M_I \to I$ be the projection. Let $\C_{M_t}$ be the constant sheaf on $M_t$ with stalk $\C$. We have then
\[ SS(\C_{M_t}) = \{ (x, t; 0, \tau) \in T^* M_I \}, \quad \sinf(\C_{M_t}) =  \{ (x, t; 0, \pm 1) \in S^* M_I \cong \tinf M\}. \]

We give another definition of isotopy of Legendrian and sheaves,  equivalent to the one given in the introduction for the case $C=S^*M$. 
\bd
Let $M$ be a smooth manifold, $I$ an open interval of $\R$. 
\begin{enumerate}
\item
An {\bf isotopy of Legendrians} over $I$ is a Legendrian $\Lambda^\infty_I\subset T^\infty (M \times I)$ such that 
\[ \linf \cap \sinf(\C_{M_t}) = \emptyset, \tforall t \in I. \] 
For any $t \in I$, we define the {\bf restriction of $\linf_I$ at $t$} as the Legendrian $\linf_t$ for the conical Lagrangian $\Lambda_t$, 
\[ \Lambda_t = \{(x,\xi) \in T^*M \mid \exists (x,t; \xi, \tau) \in \Lambda_I\}. \]

\item
An {\bf isotopy of sheaves} is a sheaf $F_I \in Sh(M \times I)$, such that 
\[ \sinf(F_I) \cap \sinf(\C_{M_t}) = \emptyset, \tforall t \in I. \] 
For any $t \in I$, we define {\bf restriction of $F_I$ at $t$} as
\[ F_t: = F_I|_{M_t} \in Sh(M). \]

\item
Two isotopies of sheaves $F_I, G_I \in Sh(M \times I)$ are {\bf non-characteristic} if 
\[ \sinf(F_t) \cap \sinf(G_t) = \emptyset, \tforall t \in I. \] 
\end{enumerate}
\ed

Some easy to check properties are in order.
\bp
(1) If $F_I$ is an isotopy of sheaf, $\linf_I = SS^\infty(F_I)$, then 
\[ \linf_t = \sinf(F_t). \]
(2) If $F_I$ is an isotopy of sheaf, $\pi_I: M_I \to I$, then $(\pi_I)_* F_I$ is a local system on $I$. 
\ep

\ss{Invariance of morphism under non-characteristic isotopy}
We use the same notations for $M_I=M \times I, M_t, \C_{M_t}, \cdots$ as in the previous subsection. 

\bl
\label{lm:nd}
Let $F \in Sh(M)$. Let $\varphi: M \to \R$ be a $C^1$ function, such that $d\varphi(x)\neq 0$ for $x \in \varphi^{-1}([0,1])$. \\
(1) For $s \in (0,1)$, let $U_s = \{x: \varphi(x) < s \}$, and let $U_1 = \cup_{s} U_s$. If 
\[ SS^\infty(\C_{U_s}) \cap SS^\infty(F) = \emptyset , \; \forall\, 0<s<1, \]
then 
\[ \Hom(\C_{U_1}, F) \congto \Hom(\C_{U_s}, F),  \; \forall\, 0<s<1. \]
(2) For $s \in (0,1)$, let $Z_s = \{x: \varphi(x) \leq s \}$, and let $Z_0 = \cap_{s} Z_s$. 
If 
\[ SS^\infty(\C_{Z_s}) \cap SS^\infty(F) = \emptyset , \; \forall\, 0<s<1, \] 
then 
\[ \Hom(\C_{Z_s}, F) \congto \Hom(\C_{Z_0}, F),\; \forall\, 0<s<1. \]
\el
\bpf
(1)  is a special case in \cite[Prop 1.8]{GKS}. (2)  follows from (1) and 
\[ 0 \to \C_{M \RM \Z_s} \to \C_M \to \C_{Z_s} \to 0. \]
\epf

The following lemma is also often used.
\bl[Petrowsky theorem for sheaves, Corollary 4.6 \cite{S}] \label{l:pet}
Let $F, G \in Sh(M)$. If $SS^\infty(F) \cap SS^\infty(G) = \emptyset$, then the natural morphism 
\[ \uhom(F, \C_M) \ot G \to \uhom(F, G) \]
is an isomorphism. 
\el
\bc
If $F_I$ be an isotopy of sheaves, then 
\[ \uhom(\C_{M_t}, F_I) \cong \C_{M_t}[-1] \ot F_I \]
\ec

\bp\label{p:nc}
Let $G_I$ and $F_I$ be non-characteristic isotopy of sheaves, then
$\uhom(F_I, G_I)$ is an isotopy of sheaves. 
In particular, 
\[ \Hom(F_t, G_t) \cong \Hom(F_s, G_s) \quad \text{ for all } t, s \in I \]
\ep
\bpf
$G_I$ and $F_I$ being non-characteristic implies $\sinf(G_I) \cap \sinf(F_I) = \emptyset$, hence we can bound singular support of the hom sheaf as \cite{KS}
\[ SS(\uhom(F_I, G_I)) \subset SS(G_I) + SS(F_I)^a. \]
Again, using $G_I$ and $F_I$ being non-characteristic, we have 
\[ SS^\infty (\uhom(F_I, G_I)) \cap SS^\infty(\C_{M_t}) = \emptyset \quad \text{ for all } t, s \in I. \]
Hence $\uhom(F_I, G_I))$ is an isotopy of sheaves. For the second statement, we have
\bee
&&  \Hom(F_t, G_t) \notag \\
&=& \Hom({j_t}^* F_I, {j_t}^* G_I) \cong \Hom(F_I, {j_t}_* j_t^* G_I) \cong \Hom(F_I, \C_{M_t} \ot G_I)   \notag\\
& \cong&   \Hom(F_I, \uhom(\C_{M_t}, G_I)[1]) \cong \Hom(\C_{M_t}, \uhom(F_I, G_I))[1]  \notag\\
&\cong& \Hom(\C_t, \pi_{I*}  \uhom(F_I, G_I))[1] \cong [\pi_{I*}  \uhom(F_I, G_I)]_t \label{eqhom}
\eee
then the result follows since $\pi_{I*}(\uhom(F_I, G_I))$ is a local system.  
\epf

\ss{Invariance of Morphism under Reeb Perturbation} \label{reeb}
Sometimes we want to vary $G, F$ while preserving $\Hom(F,G)$, but $\sinf(G) \cap \sinf(F) \neq \emptyset$, e.g. $F=G$. Here we borrow an idea from infinitesimally wrapped Fukaya-category \cite{NZ}, that to compute $\Hom_{Fuk}( L_1, L_2)$  one need to do perturbation to separate $L_1, L_2$ at infinity, one can perturb $L_2 \leadsto R^t L_2$ or $L_1 \leadsto R^{-t} L_1$ where $R^t$ is Reeb flow
\footnote{Note that in (partially) wrapped Fukaya category, one wraps $L_1$ positively (or $L_2$ negatively) in Reeb direction. This difference in sign is due to an opposite sign convention for $\omega$. Hence Reeb flow here should be termed 'geodesic flow' to be precise.}
 for positive small time $t$, small enough so that no new intersections are created between $L_1, L_2$ at infinity.\footnote{We thank P. Schapira and S. Guillermou for discussion about positive Reeb perturbation on sheaves.} 

Fix a Riemannian metric $g$ on $M$, and identify $S^*M$ with $T^\infty M$, so that Reeb flow $R^t$ is the unit speed geodesic flow. Let $r_{inj}(M, g)$ be the injective radius of $(M,g)$.
Let $\h R^t$ be the GKS quantization of $R^t$. The remaining part of this subsection will be devoted to prove the following Proposition. 
\bp\label{p:nc2}
Let $\linf \subset T^\infty M$ be a Legendrian, and $0 < \epsilon<r_{inj}(M, g)$ be small enough such that 
\[ \linf \cap R^t \linf = \emptyset , \quad \forall \; 0 < |t| < \epsilon. \] 
(1) For any $F \in Sh(M, \Lambda), 0\leq t<\epsilon$, we have canonical morphism
\[ F \to \h R^t F.\]
(2) For any $F,G \in Sh(M, \Lambda), 0 \leq t<\epsilon$, we have canonical quasi-isomorphisms
\[ \Hom(F, G) \congto \Hom(F, \h R^t G),
\quad  \Hom(F, G) \congto \Hom(\h R^{-t}F, G) \]
\ep
\bpf
For any $0\leq t < \epsilon$, define
\[ K_t=\C_{\{(x,y) | d_g(x,y) \leq t\}} \in Sh(M \times M). \]
Then from \cite{GKS}, we have 
\[ \h R^t F = \pi_{1*} \uhom(K_t, \pi_2^! F), \]
and 
\[ \h R^{-t} F = \pi_{2!} \uhom(K_t \ot \pi_1^* F), \]
where $\pi_1$ and $\pi_2$ are the projection from $M \times M$ to the first and second factor, and $\uhom$ is the (dg derived) sheaf-hom. 
From the canonical restriction morphism $K_t \to K_0=\C_\Delta$, where $\Delta \subset M \times M$ is the diagonal subset, we have 
\[ F = \pi_{1*} \uhom(K_0, \pi_2^! F) \to \pi_{1*} \uhom(K_t, \pi_2^! F) = \h R^t F. \]

For the second statement, we first prove the following   lemma. 
\bl\label{lm:Kt}
\be SS^\infty (K_t) \cap SS^\infty(\uhom(\pi_1^*F, \pi_2^! G)) = \emptyset, \quad \forall 0 < t < \epsilon \label{eq:ck}. \ee
\el
\bpf
Assuming the intersection is non-empty and contains $(x_1, x_2; p_1, p_2)$ in its cone. 
Since $(x_1, x_2; p_1, p_2) \in \R_{>0} \cdot SS^\infty(K_t)$, we have 
\[  d_g(x_1, x_2) = t. \] 
Using the boundary defining inequality $d(x_1, x_2) \leq t$ for $K_t$, we found its inward conormal at point $(x_1, x_2)$  is given by 
\[ (p_1, p_2) \in \R_{>0} \cdot (-\pa_{x_1}d_g(x_1, x_2), -\pa_{x_2} d_g(x_1, x_2)),\]
In particular, since $0 < d_g(x_1, x_2) = t <\epsilon <r_{inj}(M,g)$, $x_1, x_2$ are conjugate pairs, hence 
from the geometry of geodesic flow, we have 
\be R^t(x_1, p_1) = (x_2, -p_2), \quad
 R^t( x_2, p_2)= (x_1, -p_1), \quad p_1, p_2 \neq 0 \label{eA} \ee
   
On the other hand, since $(x_1, x_2; p_1, p_2) \in  \R_{>0} \cdot  SS^\infty(\uhom(\pi_1^*F, \pi_2^! G))$, and since $p_1, p_2 \neq 0$, we have 
\be (x_1, -p_1) \in \R_{>0} \cdot SS^\infty (F), \quad (x_2, p_2)\in \R_{>0}\cdot SS^\infty(G). \label{eB}\ee
Hence, combining \eqref{eA} and \eqref{eB}, we have 
\[ [(x_1, -p_1)] \in  R^t( SS^\infty(G)) \cap SS^\infty(F) \subset R^t \linf \cap \linf \]
This contradicts with the condition on $\epsilon$, hence finishes the proof of the Lemma. 
\epf

Now we come back to the proof of the main proposition. We have
\bea
\Hom(F, G) &\cong& \Gamma(M, \uhom(F,G)) \\
&\cong& \Gamma(M \times M, \uhom(\C_\Delta, \uhom(\pi_1^*F,\pi_2^!G)) \\
&\congto& \Gamma(M \times M, \uhom(K_t, \uhom(\pi_1^*F,\pi_2^!G)) \\
&\cong& \Gamma(M \times M, \uhom(\pi_1^*F, \uhom(K_t,\pi_2^!G)) \\
&\cong& \Gamma(M, \uhom(F, \pi_{1*}\uhom(K_t,\pi_2^!G)) \\
&\cong& \Hom(F, \h R^t G).
\eea
where in the third step when we replace $\C_\Delta$ by $K_t$, we used the canonical morphism $K_t \to \C_\Delta$, and used  Lemma \ref{lm:Kt} and Lemma \ref{lm:nd}(2) to show it is an quasi-isomorphism. 
\epf

We will use the following purely sheaf-theoretical statement later to study family of GKS quantization. 

\bp
Let $I=(0,1)$, and $K_I \in Sh(M \times M \times I)$ be an isotopy of sheaves, such that $K_t = \C_{\Delta_t}$ for some  closed subsets $\{\Delta_t\}_{0<t<1}$ satisfying
\[  \Delta_t \subset \Delta_s, \quad \forall 0<t<s<1,  \z{ and } \bigcap_{t \in I} \Delta_t = \Delta_M = \{(x,x): x \in M\}  \]
Let $F, G \in Sh(M, \Lambda)$, and $\uhom(\pi_1^*F,\pi_2^!G) \in Sh(M \times M)$ be the hom-sheaf. Assume 
\[ \sinf(K_t) \cap \sinf(\uhom(\pi_1^*F,\pi_2^!G)) = \emptyset,\quad \forall t \in I \]
then 
\[ \Hom(F,G) \cong \Hom(F, K_t^! G) \cong \Hom(K_{t!} F, G) ,\quad \forall t \in I \]
where $K_t^!, K_{t!}$ are defined in \eqref{K:shriek}. 
\ep

Its proof is exactly as in Proposition \ref{p:nc2} (2), where the condition provided in Lemma \ref{lm:Kt} is put into the hypothesis, hence we do not repeat here.

\ss{Limit of Contact Isotopy} \label{limcon}
Here we consider the compactification of $I=(0,1)$  at $0$ to $[0,1)$. Let $\Lambda_I^\infty, \linf_t$ be as before.

Denote the inclusions as 
\[  (0,1) \xhookrightarrow{j_I} [0,1) \xhookleftarrow{j_0} \{0\}. \]
\bp \label{limit}
Let $F_I$ be an isotopy of sheaves, and $\Lambda^\infty_I = SS^\infty(F_I)$. Suppose the family ${(\Lambda_t^\infty, t)} \subset T^\infty M \times (0,1)$ has a closure in $T^\infty M \times [0,1)$ whose intersection with $T^\infty M \times \{0\}$ is a Legendrian  $\Lambda_0^\infty$.  Then the sheaf
\be F_0: = (j_0)^*(j_I)_* F_I. \ee
is a constructible sheaf with $\sinf(F_0) \subset \linf_0$. 
\ep
\bpf
Suppose $(x, \xi) \notin \La_0$, with $\xi \neq 0$. We build test function $f$ in a small coordinate ball $B$ around $x$, that $f(x)=0, df(x)=\xi$. We then want to show the following 
\be\label{myiso}  F_I((B \cap \{f<\epsilon\}) \times (0, \delta)) \congto F_I((B \cap \{f<-\epsilon\}) \times (0, \delta)) \ee
for small enough $\epsilon, \delta$ and $B$. Since the limit of $\linf_t$ does not contain $[(x, \xi)]$, hence for $0<t<t_0\ll1$, we have an open conic neighborhood $\Omega \subset \dot T^*M$ of $(x,\xi) \in \dot \La_0$, such that $\La_t \cap \Omega = \emptyset$. In particular, we have 
\[ (\Omega \times T^*(0,t_0)) \cap \Lambda_I = \emptyset \]
Thus, we may choose $\epsilon, \delta$ and $B$ small enough, that the retraction $(B \cap \{f<\epsilon\}) \times (0, \delta)$ to $(B \cap \{f<-\epsilon\}) \times (0, \delta)$ is non-characteristic, hence \eqref{myiso} is an quasi-isomorphism.
\epf
\brem
We thank E. Zaslow for suggesting this condition on the family. For general behavior of how singular support of sheaves behave under pushforward or pullback of constructible sheaf, we refer the reader to \cite[Chapter 5,6]{KS}.
\erem

\section{Existence and Uniqueness of Extension} 
\bt \label{t:ext}
Let $\{\linf_t\}_{t \in I}$  be an isotopy of Legendrian in $\tinf M$, such that $\linf_t$ is constant for $t$ outside a closed interval $[a,b] \subset I$. Assume $\{\linf_t\}$ is uniformly displaceable for some Reeb vector field $R_\alpha$ on $\tinf M$,  and there exists an isotopy of convex tubular neighborhoods $\{(U(\linf_t), \rho_t, X_t)\}$ of $\{\linf_t\}$.  Denote the inclusion of slice by 
\[ \iota_t: M_t := M \times \{t\} \; \into \;M_I := M \times I. \]
Then the restriction functor
\[ \iota_t^{*}: Sh(M_I, \linf_I) \to Sh(M_t, \linf_t) \]
is an equivalence of category for all $t \in I$. 

\et
This theorem together with Proposition \ref{p:hyp2tube} implies our main theorem in the introduction. 

In the remaining part of this section, we  will sometimes identify $\linf_t \In \tinf M$ with $\lcal_t \In S^*M$, and identify Reeb flow with geodesic flow. 
\ss{Uniqueness of Extension}
\bp \label{p:uniq}
Let $\linf_t$ be a family of Legendrian in $\tinf M$ that are uniformly displaceable with parameter $\epsilon$. 
Then, the restriction functor $\iota_t^*$ is fully-faithful for all $t$. 
\ep

\bpf
For $0 \leq s < \epsilon$, we define a family of kernels in $Sh((M_1 \times I_1) \times (M_2 \times I_2))$. 
\[ K_s := \C_{d(x_1, x_2) \leq s} \boxtimes \C_{t_1=t_2}.  \]
One can check that $K_s$ generate slice-wise geodesic flow, i.e., if $F_I \in Sh(M_I)$, and
\[   K_s^! F_I:= \pi_{1*} \uhom(K_s, \pi_2^! F_I) \]
then we have
\[ SS^\infty(  (K_s^! F_I)|_{M_t} ) = R^s \sinf(F_I|_{M_t}) \] 
where $\pi_i$ is the projection from $(M_1 \times I_1) \times (M_2 \times I_2)$ to $M_i \times I_i$, and $R^s$ is the Reeb (geodesic) flow for time $s$. 

We first prove the following claim: for any $F_I, G_I \in Sh(M_I, \linf_I)$, we have 
\[ \Hom(\C_{M \times (a,b)}, \uhom(F_I, G_I)) \z { is independent of $a<b$. } \]
Suffice to prove the case for the right end-point $b$. To use the estimate of the singular support of the hom-sheaf, we would like to perturb $G_I$ by the fiberwise Reeb flow. 

\bl \label{shrink}
For any $0<s<\epsilon$, we have
\[ \Hom(\C_{M \times \{t\}}, \uhom(F_I, G_I)) \congto \Hom(\C_{M \times \{t\}}, \uhom(F_I, K_s^!G_I)). \]
The same is true if we replace $\{t\}$ by any sub-interval, eg. $[a,b]$, $(a,b)$ of $I$. 
Furthermore, $\Hom(\C_{M \times \{t\}}, \uhom(F_I, K_s^!G_I))$ is independent of $t$. 
\el
\bpf
Unwind the definition of $K_s^!$, we have
\bea
&&\Hom(\C_{M \times \{t\}}, \uhom(F_I, K_s^!G_I)) \\
& = & \Hom(\C_{M \times \{t\}}, \uhom(F_I, \pi_{1*} \uhom(K_s, \pi_2^! G_I))) \\
&=& \Hom(\C_{M \times \{t\}}, \pi_{1*} \uhom(\pi_1^* F_I,  \uhom(K_s, \pi_2^! G_I))) \\
&=&\Hom(\pi_1^* \C_{M \times \{t\}}, \uhom(K_s,  \uhom(\pi_1^* F_I, \pi_2^! G_I))) \\
\eea
 
We claim that 
\be \sinf (\pi_1^* \C_{M \times \{t\}}) \cap \sinf \uhom(K_s,  \uhom(\pi_1^* F_I, \pi_2^! G_I)) = \emptyset, \; \forall \; 0<s<\epsilon. \label{cl1} \ee

By the same argument as Proposition \ref{p:nc2} and Lemma \ref{lm:Kt}, we have 
\[ \sinf(K_s) \cap \sinf(\uhom(\pi_1^* F_I, \pi_2^! G_I)) = \emptyset,  \; \forall \; 0<s<\epsilon. \]
Hence
\be SS(\uhom(K_s,  \uhom(\pi_1^* F_I, \pi_2^! G_I))) \subset  (-\La_I, \La_I) - SS(K_s). \label {hint1}\ee
On the other hand
\bee \sinf(\pi_1^* \C_{M \times \{t\}}) = \{ [(x_1, t_1; \xi_1, \tau_1), (x_2, t_2; \xi_2, \tau_2)]: \notag \\
\xi_1=\xi_2=0, \tau_2=0, (t_1, \tau_1) = (t,\pm 1)\} \label{hint2}
\eee
If \eqref{cl1} is false, and contains a non-empty intersection point, then at the intersection we have
\[ \tau_1 + \tau_2 = \pm 1 \neq 0, \xi_1=\xi_2=0 \]
from \eqref{hint2}. From \eqref{hint1}, suppose that we  
\[ (x'_1, t'_1; \xi'_1, \tau'_1), (x'_2, t'_2; \xi'_2, \tau'_2)\in (-\La_I, \La_I) \]
where $\xi'_i = 0$ imply $\tau'_i=0$ for $i=1,2$ respectively, and 
\[  (x'_1, t'_1; \xi''_1, \tau''_1), (x'_2, t'_2; \xi''_2, \tau''_2) \in -SS(K_s) \]
where $t'_1=t'_2$, $\tau''_1+\tau''_2=0$, and if $\xi''_1 \neq 0$ iff $\xi'_2 \neq 0$ and if true implies  $d(x'_1, x'_2)=s$. 
Since $K_s|_t \circ \La_t$ is disjoint from $\La_t$ away from zero-section, hence there is no non-trivial solution to 
\[ \xi'_1 + \xi''_1 = 0, \; \xi'_2 + \xi''_2=0 \]
ie. each summand in each equation vanishes. That implies $\tau'_i=0$. Then $\tau''_1 + \tau''_2=0$ contracdicts with $\tau_1 + \tau_2 \neq 0$. Hence we proved the Lemma. 

From this claim, and
\bea \Hom(\pi_1^* \C_{M \times \{t\}}, \uhom(K_s,  \uhom(\pi_1^* F_I, \pi_2^! G_I))) 
\\
\cong \Hom(\pi_1^* \C_{M \times \{t\}} \ot K_s,  \uhom(\pi_1^* F_I, \pi_2^! G_I))) \eea
we may apply Lemma \ref{lm:nd} (2) on shrinking closed set, to get
\[ \Hom(\pi_1^* \C_{M \times \{t\}} \ot K_0,  \uhom(\pi_1^* F_I, \pi_2^! G_I))) \cong \Hom(\pi_1^* \C_{M \times \{t\}} \ot K_0,  \uhom(\pi_1^* F_I, \pi_2^! G_I))) \]
for all $0<s<\epsilon$. This proves the first statement of the Lemma. 

The final statement of the Lemma follows from \eqref{cl1}, then we may apply Proposition \ref{p:nc} where the first slot is $\C_{M \times \{t\}}$ and the second slot in hom $\uhom(F_I, K_s^!G_I)$ is taken as a constant isotopy of sheaf with any fixed $0<s<\epsilon$. The case for sub-interval can be proved similarly, and we omit the details. 
\epf
Now, we finish to prove the proposition.  By Lemma \ref{shrink}, 
\[ \Hom(\C_{M \times (a, b )}, \uhom(F_I, G_I)) \]
is independent of $(a,b)$, hence we may shrink from $(0,1)$ to an arbitrary small neighborhood of $t$. Then we have
\bea \Hom(F_I, G_I)) \cong [\pi_{I*} (\uhom(F_I, G_I)]_t \cong [\pi_{I*} (\uhom(F_I, K_s^! G_I)]_t \\
\cong \Hom(\iota_t^*F_I,\iota_t^* K_s^! G_I)   \cong  \Hom(F_t,R^s G_t) \cong  \Hom( F_t,  G_t)
\eea
where $0<s<\epsilon$, and we used small Reeb perturbation to make $F_I, K_s^! G_I$ non-charactersitic isotopy of sheaves, then apply \eqref{eqhom} in  Proposition \ref{p:nc}. 
\epf

\bp
Let $\{\linf_t\}$ be a family of Legendrian in $\tinf M$ that are uniformly displaceable with parameter $\epsilon$. For a given $t$, let $F_t \in Sh(M, \linf_t)$. Suppose we have $F'_I$ and $F''_I$ in $Sh(M_I, \linf_I)$ and isomorhpism
\[ f: F'_I|_t \congto F_t, \quad  g: F''_I|_t \congto F_t, \]
then there exist canonical isomorphism 
\[ \Phi: F'_I \to F''_I \]
such that $\Phi|_{t} = g^{-1} \circ f: F'_I|_t  \to F''_I|_t$. 
\ep
\bpf
By Proposition \ref{p:uniq}, we have  $\Hom(F'_I|_t ,  F''_I|_t) \cong \Hom(F'_I, F''_I)$. Thus, 
\[ g^{-1} \circ f \in \Hom(F'_I|_t ,  F''_I|_t) \mapsto \Phi \in  \Hom(F'_I, F''_I). \] 
Similarly 
\[ f \circ g^{-1} \in \Hom(F''_I|_t ,  F'_I|_t) \mapsto \Psi \in \Hom(F''_I, F'_I). \]
Hence we have
\[ \Phi_t \circ \Psi_t \cong id_{F''_t} \in \Hom(F''_I|_t ,  F''_I|_t) \mapsto \Phi  \circ \Psi  \cong id_{F''_I} \in \Hom(F''_I ,  F''_I). \] 
and similarly $ \Psi  \circ \Phi  \cong id_{F'_I}$. 
\epf

\ss{Existence of Local Extension}

\bp\label{p:localext}
Let $\{\lcal_t\}$ be a family of Legendrian in $S^*M$ that admits a family of convex tubular neighborhood thickening $\{(U,\rho, X)_t\}$. Then  for any compact subset $K \subset I$, there exists $\delta>0$ such that for  any $t \in K$ and $F_t \in Sh(M, \lcal_t)$,  there exists $F_J \in Sh(M \times J, \lcal_{J})$ where $J=I \cap (t-\delta, t+\delta) \subset I$, such that $F_J|_t \congto F_t$ canonically. 
\ep
\bpf
Define $Q$ as neighborhood of diagonal in $I \times I$
\[ Q = \{(s,t) \in I \times I \mid \lcal_s \in U_{1/2} (\lcal_t), \;  \lcal_t \in U_{1/2} (\lcal_s),\} \]
Then, we may find $\delta=\delta(K, \{U_t, \rho_t\})$ small enough such that 
\[ \Delta_{K,\delta} = \bigcup_{t \in K} [t-\delta, t+\delta] \times [t-\delta, t+\delta]\]
is contained in $Q$. 

For any $t \in K$, let $J=I \cap (t-\delta, t+\delta) \subset I$. For any $s \in J$, we consider the trajectory of $\lcal_t$ under the retracting flow $-X_s$, and get an isotopy of Legendrians over $[0,\infty)$ as $X_s^{-T}(\lcal_t)$. We claim that the Gromov-Hausdorff limit of $X_s^{-T}(\lcal_t) \In S^*M$ is $\lcal_s$, since 
\[ X_s^{-T}(\lcal_t) \subset X_s^{-T}(U_{1/2}(\lcal_s)) \subset U_{e^{-c_s T}/2}(\lcal_s) \to \lcal_s, \z{ as } T \to \infty, \]
where $c_s$ is the shrinking rate $\la d \rho_s, X_s \ra > c_s \rho_s$ in the definition of $(U,\rho,X)$. 
Thus we may define the limit of the corresponding isotopy of sheaves 
\[ \Pi_s(F_t): = (j_\infty)^* (j_{[0,\infty)})_* (\h X_s^{-[0,\infty)} F_t) \]
where 
\[ \h X_s^{-[0,\infty)}: Sh(M) \to Sh(M \times [0,\infty))\] 
is the sheaf quantization of the flow $X_s^{-T}$ and
\[ j_{[0,\infty)}: \; [0,\infty) \hookrightarrow \;[0,\infty]\; \hookleftarrow \{\infty\}: \; j_{\infty} \]
are inclusion into the compactification, and we also abuse notation to denote $id_M \times j$ as $j$. By Proposition
\ref{limit}, we have 
\[ \sinf(\Pi_s(F_t)) \subset \lcal_s. \] 

We claim that the collection of sheaves $\{\Pi_s(F_t)\}_{s \in J}$ assemble into an isotopy of sheaf, $\Pi_J(F_t) \in Sh(M_J, \lcal_J)$. Indeed since the contact flow $X_s$ varies smoothly with parameter $s$, we have a (tensor) kernel for the family
\[ K_J \in Sh(M_J \times (M_J \times [0,\infty))), \]
such that we have 
\[ \Pi_J(F_t): = (id_{M_J} \times j_\infty)^* (id_{M_J} \times j_{[0,\infty)})_*( (K_J)_! (F_t \boxtimes \C_J)). \]
Thus we get the extension sheaf $\Pi_J(F_t)$, and one can check $\Pi_J(F_t)|_t \cong F_t$ since the retraction flow $X_t^{-T}$ preserves the Legendrian $\lcal_t$. 
\epf

\ss{Proof of Theorem \ref{t:ext}}
Let $K = [a,b] \subset (0,1)$, and we apply Proposition \ref{p:localext} to get the positive constant $\delta>0$, such that for any $t \in K$, we may extend a sheaf $F_t \in Sh(M, \lcal_t)$ to a neighborhood $B_\delta(t)=(t-\delta, t+\delta)$, compatible with the Legendrian condition $\lcal_I$ restricted on the interval. We may take a finite set of points 
\[ A=\{t_n \in I \mid t_n = t + (n/2) \delta/2, n \in \Z, t_n \in [a,b] \} \]
and extend the sheaf $F_t$ from $M \times \{t\}$ inductively to $M \times B_{ (1+n/2)\delta}(t)$ for $n=0,1,2,\cdots$, using existence of local extension and uniqueness of extension. Finally, since the isotopy is constant outside $[a,b]$, we may trivially extend from $[a,b]$ to $(0,1)$. This finishes the proof of Theorem \ref{t:ext}. 
%
%
%
%
%
%
%



\begin{thebibliography}{HHHH}

\bibitem[Av]{Av} Russell Avdek, Liouville hypersurfaces and connect sum cobordisms, arXiv:1204.3145.

\bibitem[CE]{CE} K. Cieliebak and Y. Eliashberg. From Stein to Weinstein and Back: Symplectic Geometry and Affine Complex Manifolds. AMS Colloquium Publications, vol. 59 (2012)

\bibitem[Eli]{Eli} Y. Eliashberg. Weinstein manifolds revisited. arXiv:1707.03442 




\bibitem[Ge]{Ge} H. Geiges,  An Introduction to Contact Topology. 

\bibitem[Gi]{Gi} E. Giroux, Convexite en topologie de contact, Comment. Math. Helv. 66 (1991), 637–677.

\bibitem[GKS]{GKS} S. Guillermou, M. Kashiwara, P. Schapira. Sheaf quantization of Hamiltonian isotopies and applications to nondisplaceability problems.  Duke Math. J. Volume 161, Number 2 (2012), 201-245.

\bibitem[GPS]{GPS} S. Ganatra, J. Pardon, V. Shende. Covariantly functorial wrapped Floer theory on Liouville sectors. arXiv:1706.03152


\bibitem[KS] {KS} M. Kashiwara, Pierre Schapira, Sheaves on Manifolds



 
 

\bibitem[N1]{N1} D. Nadler. Microlocal branes are constructible sheaves, Selecta Math. 15 (2009), no. 4, 563--619.


\bibitem[N2]{N2} D. Nadler. Arboreal Singularities.  Geometry \& Topology 21 (2017) 1231 --1274


\bibitem[N3]{N3} D. Nadler. Non-characteristic expansions of Legendrian singularities. arXiv:1507.01513

\bibitem[N4]{N4} D. Nadler. Wrapped microlocal sheaves on pairs of pants. arXiv:1604.00114

\bibitem[NZ]{NZ} D. Nadler, E. Zaslow. Constructible Sheaves and the Fukaya Category. J. Amer. Math. Soc. 22 (2009), 233-286 


\bibitem[S]{S} Pierre Schapira. A short review on microlocal sheaf theory. \href{https://webusers.imj-prg.fr/~pierre.schapira/lectnotes/MuShv.pdf}{link}. 

\bibitem[Sh]{Sh} V. Shende. Microlocal category for Weinstein manifolds via h-principle.  arXiv:1707.07663 

\bibitem[Syl]{Syl} Z. Sylvan. On partially wrapped Fukaya categories. arXiv:1604.02540
\end{thebibliography}
\end{document}